\documentclass[12pt]{article}
\title{ Stability of the line soliton of the Kadomtsev--Petviashvili-I equation with the critical traveling speed 
\amssubj{35B35, 37K40, 35Q53}
}
\author{Yohei Yamazaki \footnote{{\it E-mail addresses:} yohei.yamazaki@u-cergy.fr} \\ {\footnotesize University of Cergy-Pontoise,
 } \\  {\footnotesize  UMR CNRS 8088, Cergy-Pontoise, F-95000, France}\\ {\footnotesize and } \\ {\footnotesize Osaka City University Advanced Mathematical Institute }\\ {\footnotesize 3-3-138 Sugimoto, Sumiyoshi-ku Osaka 558-8585 Japan }}
\usepackage{amssymb}
\usepackage{mathrsfs}
\usepackage{amsthm}
\usepackage{amsmath}
\usepackage{titlefoot}
\usepackage{bm}
\def\pdfliteral #1 {}

\numberwithin{equation}{section}

\newtheorem{theorem}{Theorem}[section]

\newtheorem{lemma}[theorem]{Lemma}

\theoremstyle{definition}

\newtheorem{remark}[theorem]{Remark}

\renewcommand{\eqref}[1]{(\ref{#1})}

\renewcommand{\bigskip}{\vspace{0.3cm}}

\newcommand{\R}{{\mathbb R}}

\newcommand{\Z}{{\mathbb Z}}
\newcommand{\T}{{\mathbb T}}

\newcommand{\norm}[1]{{\left \lVert #1 \right \rVert}}

\newcommand{\rbr}[1]{\left( #1 \right)}
\newcommand{\tbr}[1]{\langle #1 \rangle}

\newcommand{\Tbr}[1]{\left\langle #1 \right\rangle}
\newcommand{\RT}{\mathbb {R} \times \mathbb {T}}

\date{}
\begin{document}

\maketitle

{\small keywords and phrases: Nonlinear stability,  line soliton, Kadomtsev--Petviashvili-I equation}

\begin{abstract}
We consider the orbital stability of solitons of the Kadomtsev--Petviashvili-I equation in $\R \times (\R/2\pi\Z)$ which is one of a high dimensional generalization of the Korteweg--de Vries equation.
In \cite{TBB}, Benjamin showed that the Korteweg--de Vries equation possesses the stable one soliton.
We regard the one soliton of the Korteweg--de Vries equation as a line soliton of  the Kadomtsev--Petviashvili-I equation.
Zakharov \cite{VEZ} and Rousset--Tzvetkov \cite{R T 0} proved the orbital instability of the line solitons of the Kadomtsev--Petviashvili-I equation on $\R^2$.
The orbital instability of the line solitons of the Kadomtsev--Petviashvili-I equation on $\R \times (\R/2\pi\Z)$ with the traveling speed $c>4/\sqrt{3}$ was proved by Rousset--Tzvetkov \cite{R T 1} and the orbital stability of the line solitons with the traveling speed $0<c<4/\sqrt{3}$ was showed in \cite{R T 3}.
In this paper, we prove the orbital stability of the line soliton of the Kadomtsev--Petviashvili-I equation on $\R \times (\R/2\pi\Z)$ with the critical speed $c=4/\sqrt{3}$ and the Zaitsev solitons near the line soliton.
Since the linearized operator around the line soliton with the traveling speed $4/\sqrt{3}$ is degenerate, we can not apply the argument in \cite{R T 1,R T 3}.
To prove the stability of the line soliton, we investigate the branch of the Zaitsev solitons.
\end{abstract}

\section{Introduction}

% We consider the stability of the line soliton of KP-I equation with the critical traveling speed. 
We consider the Kadomtsev--Petviashvili-I equation:
\begin{align}\label{KP-I-eq}
(u_{t}+ u_{xxx} + uu_x)_x -u_{yy}=0 \quad (t,x,y)\in \R \times \RT,
\end{align}
where $\T=\R/2\pi \Z$.
In \cite{K P}, Kadomtsev and Petviashvili introduced the equation \eqref{KP-I-eq} as a model equation for  the propagation of long waves weakly modulated in the transverse direction.
The local and global well-posedness of the KP-I equation on $\R^2$ was showed by \cite{G P W,I K T,CK,M S T 0,SU}.
The Cauchy problem of the KP-I equation for initial datum which are localized perturbations of a non-localized traveling wave solution was studied by \cite{F P,M S T}.
In \cite{I K}, Ionescu and Kenig proved the global well-posedness of the KP-I equation on $\RT$ for initial data in the second energy space $ \mathcal{Z}^2$, where 
\[\mathcal{Z}^s=\{ u  : \norm{u}_{\mathcal{Z}^s}<\infty\}\]
and
\[\norm{u}_{\mathcal{Z}^s}=\norm{(1+|\xi|+|k\xi^{-1}|)^s\hat{u}(\xi,k)}_{L^2(\R_{\xi}\times\Z_{k})}.\]
The conservation laws
\[M(u)=\int_{\RT}|u|^2dxdy\]
and
\[E(u)=\int_{\RT}\Bigl( |\partial_x u|^2 + |\partial_x^{-1}\partial_y u|^2 - \frac{1}{3}u^3 \Bigr)dxdy\] of \eqref{KP-I-eq} identify the space $\mathcal{Z}^1$ as the energy space of \eqref{KP-I-eq}.

Solitons are nontrivial traveling wave solutions of the KP-I equation with the form $u(t,x,y)=Q(x-ct,y)$. 
The function $Q(x-ct,y)$ is soliton of the KP-I equation on $\RT$ if and only if  the function $Q$ is a nontrivial solution to the stationary equation
\begin{align}\label{S-eq}
-u_{xx} + \partial_x^{-2} u_{yy} +c u -\frac{1}{2} u^2 =0, \quad (x,y) \in \RT.
\end{align}
% We define the action 
% \[S_c(u)=E(u)+cM(u).\]
In \cite{M Z B I M}, Manakov et. al showed that the KP-I equation possesses the lump solitons 
\[\phi_c(x-ct,y)=\frac{8c(1-\frac{c}{3}(x-ct)^2+\frac{c^2}{3}y^2)}{(1+\frac{c}{3}(x-ct)^2+\frac{c^2}{3}y^2)^2}.\]
In \cite{dB S 1,dB S 2,dB S 3}, de Bouard and Saut proved the existence and the stability of the ground states of the stationary equation of the generalized KP-I equation on $\R^2$.

The KP-I equation is one of a high dimensional generalization of the Korteweg--de Vries equation:
\[u_t+u_{xxx}+uu_x=0, \quad (t,x) \in \R\times \R.\]
The KdV equation possesses the one solitons:
\[ Q_c(x-ct)=3c\cosh^{-2}\Bigl(\frac{\sqrt{c}(x-ct)}{2} \Bigr).\]
The orbital stability of the one solitons of the KdV equation was showed by Benjamin \cite{TBB}.
The asymptotic stability of the one solitons of the KdV equation on weighted space was proved by Pego and Weinstein \cite{P W 2} and Mizumachi \cite{TM1}.
In \cite{M M 1,M M 2}, Martel and Merle showed the asymptotic stability of the one solitons of the generalized KdV equation for perturbations on the energy space $H^1(\R)$ by establishing the Liouville theorem.
In \cite{K S,M S T,AAZa}, it is showed that the stationary equation \eqref{S-eq} has the Zaitsev solitons 
\[Z(a)(x,y)=\frac{12(1-a^2)\Bigl(1-a\sqrt{2-a^2}\cosh \Bigl( \frac{\sqrt{1-a^2}}{3^{1/4}}x \Bigr) \cos y \Bigr)}{\sqrt{3}\Bigl( \cosh\Bigl(\frac{\sqrt{1-a^2}}{3^{1/4}} x\Bigr) -a\sqrt{2-a^2} \cos y \Bigr)^2}\]
with the traveling speed 
\[c(a)=\frac{4-2a^2+a^4}{\sqrt{3}(1-a^2)}.\]
Then, $Z(0)=Q_{4/\sqrt{3}}$.
% Therefore, Zaitsev solitons comprise the branch of soliton which is connected with the branch of line solitons at $c=4/\sqrt{3}$.

In this paper, we regard the one solitons of the KdV equation as line solions of the KP-I equation on $\RT$.
By using the Lax pair structure of the KP-I equation, Zakharov showed the instability of the line soliton as the KP-I equation on $\R^2$.
The spectral stability of the line solitons as the KP equation was obtained by Alexander--Pego--Sachs \cite{A P S}.
Using the method developed by Grenier \cite{EG}, Rousset--Tzvetkov \cite{R T 0} proved the orbital instability of the line solitons as the generalized KP-I equation for localized perturbations on $\R^2$. 
In \cite{R T 1,R T 3}, Rousset and Tzvetkov showed that the line soliton $Q_c$ is orbitally stable as the KP-I equation on $\RT$ for $0<c<\frac{4}{\sqrt{3}}$ and orbitally unstable for $c> \frac{4}{\sqrt{3}}$.
To show the instability of line solitons of the generalized KP-I equation, Rousset and Tzvetkov proved the existence of a strongly stable manifolds of line solitons associated to the most negative eigenvalue of the linearized operator around the line soliton in \cite{R T 3}.
The spectral instability of periodic traveling wave solutions of the generalized KdV equation as the generalized KP equation was showed by Johnson and Zumbrun \cite{J Z}.
The stability of the line soliton as the KP-II equation was confirmed the heuristic analysis by Kadomtsev and Petviashvili \cite{K P}.
The stability of the line solitons as the KP-II equation for decaying perturbations was showed by Villarroel and Ablowitz \cite{V A}. 
Mizumachi and Tzvetkov proved the orbital stability and the asymptotic stability of the line solitons as the KP-II equation in $L^2(\RT)$ by applying the B\"acklund transformation.
The asymptotic stability of the line solitons as the KP-II equation on $\R^2$ has proved by Mizumachi \cite{TM2,TM3}.

Now let us introduce our result.
\begin{theorem}\label{thm-main}
There exists $a_0>0$ such that for $0\leq a <a_0$ the soliton $Z(a)(x-c(a)t,y)$ is orbitally stable as a solution of the KP-I  equation $\eqref{KP-I-eq}$.
More precisely, for every $\varepsilon >0$, there exists $\delta>0$ such that for $u_0 \in \mathcal{Z}^2$ with 
\[\norm{u_0-Z(a)}_{\mathcal{Z}^1}<\delta, \]
the solution $u$ to the KP-I equation $\eqref{KP-I-eq}$, defined by \cite{I K} with the initial data $u_0$ satisfies 
\[\sup_{t\geq 0} \inf_{(x_0,y_0) \in \RT} \norm{u(t,\cdot,\cdot)-Z(a)(\cdot-x_0,\cdot-y_0)}_{\mathcal{Z}^1} < \varepsilon .\]
\end{theorem}
\begin{remark}
Since $Z(0)=Q_{4/\sqrt{3}}$ and $c(0)=4/\sqrt{3}$,  Theorem $\ref{thm-main}$ yields the orbital stability of the line soliton with the critical speed $4/\sqrt{3}$.
\end{remark}
The proof of the stability of the Zaitsev solitons $Z(a)$ for $0<a<a_0$ follows the Lyapunov function method in \cite{G S S 1}.
For $0<c<4/\sqrt{3}$, the orbital stability of the line solitons as the KP-I equation follows the coercive type estimate of the linearized operator around the line solitons in \cite{R T 3}.
In the case $c=4/\sqrt{3}$, the linearized operator around the line soliton $Q_{4/\sqrt{3}}$ is degenerate.
Therefore,  to prove the orbital stability of $Q_{4/\sqrt{3}}$, we can not apply the argument in \cite{G S S 1,R T 3}.
In \cite{R T 1}, Rousset--Tzvetkov showed the spectral instability of line solitons as the KP-I equation implies the orbital instability of the line solitons.
Since the line soliton $Q_{4/\sqrt{3}}$ is spectrally stable, we can not show the instability of $Q_{4/\sqrt{3}}$ by applying the arguments in \cite{G O,R T 3}.
To cover the degeneracy of the linearized operator of the stationary equation \eqref{S-eq}, we apply the argument of \cite{YY2,YY4} for the stability of line solitary waves with a critical exponent.
The line solitary waves with a critical exponent is a bifurcation point where the branch of line solitary waves and the branch of two dimensional solitary wave connect.
In \cite{YY2,YY4}, to prove the stability of line solitary waves with a critical exponent, the author showed the positivity of the second order term of $L^2$-norm of the two dimensional solitary waves with respect to the bifurcation parameter.
In the case \cite{YY2,YY4}, since the second order term of $L^2$-norm of the two dimensional solitary waves is not zero, we can cover the degeneracy of the linearized operator around the line solitary wave with the critical exponent by the forth order term of a Lyapunov function.
In our case, the branch of the line solitons $Q_c$ and the branch of the Zaitsev solitons $Z(a)$ connect at the bifurcation point $Q_{4/\sqrt{3}}$.
However, the second order term of $L^2$-norm of the Zaitsev solitons $Z(a)$ around $a=0$ is zero.
To cover the degeneracy of the linearized operator of \eqref{S-eq}, we use the forth order term of $L^2$-norm of the Zaitsev solitons $Z(a)$ around $a=0$ and the sixth order term of a Lyapunov function.

Our plan of the present paper is as follows.
In Section 2, we show the increase of the $L^2$-norm of the Zaitsev solitons with respect to the bifurcation parameter $a$ and prove the estimate of the Lyapunov function around the Zaitsev solitons and the line soliton $Q_{4/\sqrt{3}}$.
In Section 3, we prove the orbital stability of the line soliton $Q_{4/\sqrt{3}}$ with the critical traveling speed. 
In Section 4, we show the orbital stability of the Zaitsev solitons near the line soliton $Q_{4/\sqrt{3}}$.
\section{Preliminaries}
In this section, we show the estimate of the Lyapunov function around the Zaitsev solitons and the line soliton $Q_{4/\sqrt{3}}$.

We show the increase to the $L^2$-norm of the Zaitsev solitons with respect to the parameter $a$.
\begin{lemma}\label{lem-0}
The following hold.
\[ \partial_a\norm{Z(a)}_{L^2}^2|_{a=0}=\partial_a^2\norm{Z(a)}_{L^2}^2|_{a=0}=\partial_a^3\norm{Z(a)}_{L^2}^2|_{a=0}=0\]
and
\[\partial_a^4\norm{Z(a)}_{L^2}^2|_{a=0}=4^3 \cdot 3^{13/4}\pi  \int_{\R}f^{-4}dx,\]
where
\[\partial_a^n\norm{Z(a)}_{L^2}^2|_{a=0}=\lim_{a \to 0}\partial_a^n\norm{Z(a)}_{L^2}^2.\]
\end{lemma}
\proof
Let $\beta=a\sqrt{2-a^2}$, $f(x)=\cosh x$ and $g(y)=\cos y$.
Then, we have
\[ 1-a^2=\sqrt{1-\beta^2},\]
\begin{align}\label{eq-lem-0-1}
2k\int_{\R} f^{-2k}(x)dx=(2k+1)\int_{\R}f^{-2(k+1)}(x)dx
\end{align}
and
\begin{align*}
\norm{Z(a)}_{L^2}^2
% =& \frac{12^2 (1-\beta^2)}{3}\int_{\RT} \frac{\Bigl(1-\beta f\Bigl(\frac{(1-\beta^2)^{1/4}x}{3^{1/4}} \Bigr) g(y) \Bigr)^2dxdy}{\Bigl(f\Bigl(\frac{(1-\beta^2)^{1/4}x}{3^{1/4}} \Bigr)-\beta  g(y) \Bigr)^4}\\
=&\frac{12^2(1-\beta^2)^{3/4}}{3^{3/4}} \int_{\RT} \frac{(1-\beta f(x)g(y))^2dxdy}{(f(x)-\beta g(y))^4}.
\end{align*}
By the elementally calculation, we obtain the following equations.
\begin{align*}
\frac{3^{3/4}}{12^2}\partial_{\beta}\norm{Z(a)}_{L^2}^2
% =& -\frac{3}{2}\beta(1-\beta^2)^{-1/4}  \int_{\RT}\frac{(1-\beta fg)^2 dxdy}{(f-\beta g)^4}\\
% &+(1-\beta^2)^{3/4} \int_{\RT} \frac{-2 fg (1-\beta fg)(f-\beta g) + 4g (1-\beta fg)^2}{(f-\beta g)^5}dxdy\\
% =& -\frac{3}{2}\beta(1-\beta^2)^{-1/4}  \int_{\RT}\frac{(1-\beta fg)^2 dxdy}{(f-\beta g)^4}\\
% &+(1-\beta^2)^{3/4} \int_{\RT} \frac{-2f^2g+2\beta f^3g^2+2\beta fg^2 -2\beta^2 f^2g^3+4g-8\beta fg^2 + 4 \beta^2 f^2g^3}{(f-\beta g)^5}dxdy,
=& -\frac{3\beta}{2(1-\beta^2)^{1/4}}  \int_{\RT}\frac{(1-\beta fg)^2 dxdy}{(f-\beta g)^4}\\
&+(1-\beta^2)^{3/4} \int_{\RT}\frac{4g-2f^2g+2\beta f^3g^2-6\beta fg^2 +2\beta^2 f^2g^3}{(f-\beta g)^5}dxdy.
\end{align*}

\begin{align*}
&\frac{3^{3/4}}{12^2}\partial_{\beta}^2\norm{Z(a)}_{L^2}^2\\
% =&-\frac{3}{2}\Bigl( (1-\beta^2)^{-1/4}+\frac{\beta^2}{2}(1-\beta^2)^{-5/4} \Bigr)  \int_{\RT}\frac{(1-\beta fg)^2 dxdy}{(f-\beta g)^4}\\
% &- 3\beta (1-\beta^2)^{-1/4} \int_{\RT} \frac{4g-2f^2g+2\beta f^3g^2-6\beta fg^2 +2\beta^2 f^2g^3}{(f-\beta g)^5}dxdy\\
% &+ (1-\beta^2)^{3/4} \int_{\RT} \frac{1}{(f-\beta g)^6}\Bigl( (2f^3g^2-6fg^2+4\beta f^2 g^3)(f-\beta g)\\
% &+ 5g (4g-2f^2g+2\beta f^3g^2 -6\beta fg^2 +2\beta^2 f^2 g^3) \Bigr)dxdy\\
% =&-\frac{3}{4}(2-\beta^2)(1-\beta^2)^{-5/4}\int_{\RT}\frac{(1-\beta fg)^2 dxdy}{(f-\beta g)^4}\\
% &- 3\beta (1-\beta^2)^{-1/4} \int_{\RT} \frac{4g-2f^2g+2\beta f^3g^2-6\beta fg^2 +2\beta^2 f^2g^3}{(f-\beta g)^5}dxdy\\
% &+ (1-\beta^2)^{3/4} \int_{\RT} \frac{1}{(f-\beta g)^6}\Bigl( 2f^4g^2-6f^2g^2 +4 \beta f^3 g^3 -2\beta f^3 g^3 + 6 \beta fg^3 -4 \beta^2 f^2 g^4 \\
% &+ 20g^2 -10f^2 g^2 +10\beta f^3 g^3 -30\beta f g^3 + 10 \beta^2 f^2 g^4 \Bigr)dxdy.
=&-\frac{3(2-\beta^2)}{4(1-\beta^2)^{5/4}}\int_{\RT}\frac{(1-\beta fg)^2 dxdy}{(f-\beta g)^4}\\
&- \frac{3\beta}{ (1-\beta^2)^{1/4}} \int_{\RT} \frac{4g-2f^2g+2\beta f^3g^2-6\beta fg^2 +2\beta^2 f^2g^3}{(f-\beta g)^5}dxdy\\
&+ (1-\beta^2)^{3/4} \int_{\RT}   \frac{20g^2 + 2f^4g^2-16f^2g^2 +12 \beta f^3 g^3  -24 \beta fg^3 +6 \beta^2 f^2 g^4}{(f-\beta g)^6}dxdy.
\end{align*}
Thus, by the equation \eqref{eq-lem-1-1} we have
\begin{align*}
\frac{3^{3/4}}{12^2}\partial_{\beta}^2\norm{Z(a)}_{L^2}^2|_{\beta=0}
% =&-\frac{3}{2} \int_{\RT} \frac{1}{f^4}dxdy +\int_{\RT}\frac{20g^2+2f^4g^2-16f^2g^2}{f^6}dxdy\\
=&\pi \int_{\R} \Bigl(- \frac{3}{f^4} + \frac{20}{f^6} +\frac{2}{f^2}-\frac{16}{f^4} \Bigr)dx
=0.
\end{align*}
Similarly, we obtain that
\begin{align*}
&\frac{3^{3/4}}{12^2}\partial_{\beta}^3\norm{Z(a)}_{L^2}^2\\
% =&-\frac{3}{4}\Bigl( -2\beta(1-\beta^2)^{-5/4} +\frac{5}{2}\beta(2-\beta^2)(1-\beta^2)^{-9/4} \Bigr)\int_{\RT}\frac{(1-\beta fg)^2 dxdy}{(f-\beta g)^4}\\
% &-\frac{9}{4}(2-\beta^2)(1-\beta^2)^{-5/4}\int_{\RT}  \frac{4g-2f^2g+2\beta f^3g^2-6\beta fg^2 +2\beta^2 f^2g^3}{(f-\beta g)^5}dxdy\\
% &- \frac{9}{2} \beta (1-\beta^2)^{-1/4} \int_{\RT} \frac{1}{(f-\beta g)^6}\Bigl( 20g^2 + 2f^4g^2-16f^2g^2 +12 \beta f^3 g^3  -24 \beta fg^3 +6 \beta^2 f^2 g^4  \Bigr)dxdy\\
% &+(1-\beta^2)^{3/4}\int_{\RT} \frac{1}{(f-\beta g)^7} \Bigl( (12f^3g^3 -24 fg^3 + 12\beta f^2g^4)(f-\beta g)\\
% &+6g(20g^2+2f^4g^2-16f^2g^2+12 \beta f^3 g^3 -24\beta f g^3 + 6\beta^2 f^2 g^4) \Bigr)dxdy\\
=&-\frac{3(6\beta-\beta^3)}{8(1-\beta^2)^{9/4}} \int_{\RT}\frac{(1-\beta fg)^2 dxdy}{(f-\beta g)^4}\\
&-\frac{9(2-\beta^2)}{4(1-\beta^2)^{5/4}}\int_{\RT}  \frac{4g-2f^2g+2\beta f^3g^2-6\beta fg^2 +2\beta^2 f^2g^3}{(f-\beta g)^5}dxdy\\
&- \frac{9\beta}{2 (1-\beta^2)^{1/4}}  \int_{\RT} \frac{20g^2 + 2f^4g^2-16f^2g^2 +12 \beta f^3 g^3  -24 \beta fg^3 +6 \beta^2 f^2 g^4}{(f-\beta g)^6}dxdy\\
&+(1-\beta^2)^{3/4}\int_{\RT} \frac{120g^3 + 24f^4g^3 -120f^2g^3+72\beta f^3g^4-120\beta fg^4 +24 \beta^2 f^2 g^5}{(f-\beta g)^7}dxdy
\end{align*}
and
\begin{align*}
&\frac{3^{3/4}}{12^2}\partial_{\beta}^4\norm{Z(a)}_{L^2}^2\\
=&-\frac{9(4+12\beta^2-\beta^4 )}{16(1-\beta^2)^{13/4}} \int_{\RT}\frac{(1-\beta fg)^2 dxdy}{(f-\beta g)^4}\\
&-\frac{3(6\beta-\beta^3)}{2(1-\beta^2)^{9/4}}\int_{\RT} \frac{4g-2f^2g+2\beta f^3g^2-6\beta fg^2 +2\beta^2 f^2g^3}{(f-\beta g)^5}dxdy\\
&-\frac{9(2-\beta^2)}{2(1-\beta^2)^{5/4}} \int_{\RT} \frac{20g^2 + 2f^4g^2-16f^2g^2 +12 \beta f^3 g^3  -24 \beta fg^3 +6 \beta^2 f^2 g^4}{(f-\beta g)^6}dxdy\\
&- \frac{6 \beta}{ (1-\beta^2)^{1/4}}\int_{\RT} \frac{120g^3 + 24f^4g^3 -120f^2g^3+72\beta f^3g^4-120\beta fg^4 +24 \beta^2 f^2 g^5}{(f-\beta g)^7} dxdy\\
&+ (1-\beta^2)^{3/4}\int_{\RT} \frac{840g^4 +240f^4g^4-960f^2g^4+480\beta f^3 g^5 -720\beta fg^5 +120 \beta^2 f^2 g^6}{(f-\beta g)^8}dxdy.
\end{align*}

Since
\[\partial_a\norm{Z(a)}_{L^2}^2|_{a=0}=\partial_a^2\norm{Z(a)}_{L^2}^2|_{a=0}=\partial_a^3\norm{Z(a)}_{L^2}^2|_{a=0}=0,\]
from the equation \eqref{eq-lem-0-1} we obtain
\begin{align*}
\partial_a^4\norm{Z(a)}_{L^2}^2|_{a=0}= \Bigl(\frac{d\beta}{da}\Bigr)^4|_{a=0}\partial_{\beta}^4\norm{Z(a)}_{L^2}^2|_{\beta=0}
% =& \frac{4\cdot 12^2}{3^{3/4}} \int_{\RT} \Bigl( -\frac{9}{4f^4} -\frac{9(20g^2+2f^4g^2-16f^2g^2)}{f^6} + \frac{840g^4+240f^4g^4 -960f^2g^4}{f^8}\Bigr)dxdy\\
% =&\frac{4\cdot 12^2\pi }{3^{3/4}} \int_{\R} \Bigl(-\frac{9}{2f^4}-\frac{9(20+2f^4-16f^2)}{f^6} + \frac{630+180f^4 -720f^2}{f^8}\Bigr)dxdy\\
% =&\frac{4\cdot 12^2\pi }{3^{3/4}} \int_{\R} \Bigl(-\frac{9}{f^4}-\frac{360+36f^4-288f^2}{f^6} + \frac{1260+360f^4-1440f^2}{f^8}\Bigr)dx\\
=&\frac{2\cdot 12^2\pi }{3^{3/4}} \int_{\R} \Bigl(-\frac{36}{f^2} + \frac{639}{f^4}-\frac{1800}{f^6} + \frac{1260}{f^8}\Bigr)dx\\
=&\frac{2\cdot 12^2 \cdot 3^2\pi }{3^{3/4}} \int_{\R}f^{-4}dx.
\end{align*}
\qed

% Since the traveling speed satisfies
% \[\partial_ac(a)|_{a=0}=0 \mbox{ and } \partial_a^2c(a)|_{a=0}=\frac{4}{\sqrt{3}}>0.\]
% for small $a>0$ the Zaitsev soliton $Z(a)$ is orbitally unstable and the line soliton maybe be orbitally unstable.

Let $v_*(x)\cos y=\partial_aZ(a)|_{a=0}(x,y)$ and
\[Z(\bm{a},\gamma)(x,y)=\gamma Z(|\bm{a}|)\Bigl(\sqrt{\gamma}x,y+\theta(\bm{a})\Bigr),\]
where $(\cos \theta(\bm{a}),\sin \theta (\bm{a}))=(\frac{a_1}{|\bm{a}|}, \frac{a_2}{|\bm{a}|})$ for $\bm{a}=(a_1,a_2) \in \R^2 \setminus \{(0,0)\}$.
Then, we have
\[\partial_x^{-1}v_*(x)=\frac{12\sqrt{2}\sinh\Bigl(\frac{x}{3^{1/4}}\Bigr)}{3^{1/4}\cosh^2\Bigl(\frac{x}{3^{1/4}}\Bigr)}.\]
We define the action as
\[S_c(u)=E(u)+cM(u).\]
Let $(\cdot, \cdot)_{L^2(X)}$ be the inner product of $L^2(X)$.
% Let 
% \[\norm{u}_{\mathcal{Z}^s}=\norm{(|\xi|^s+|\xi^{-1}k|^s+1)\hat{u}(\xi,k)}_{L^2_{\xi,k}(\R \times \Z)}\]
% and
% \[\mathcal{Z}^s=\{ u \in L^2(\RT): \norm{u}_{\mathcal{Z}^s}<\infty\}.\]
The coerciveness of $S_{4/\sqrt{3}}''(Q_{4/\sqrt{3}})$ on a subspace of $\mathcal{Z}^1$ yields the following lemma.
\begin{lemma}\label{lem-5}
There exist $k_1>0$ and $\varepsilon _0>0$ such that for $a_1,a_2 \in (-\varepsilon _0,\varepsilon _0)$ and $\gamma \in (1-\varepsilon _0, 1+\varepsilon _0)$, if $w \in \mathcal{Z}^1$ satisfies
\begin{align*}
&\rbr{w,Z(\bm{a},\gamma)}_{L^2(\RT)}=\rbr{w,\partial_x Z(\bm{a},\gamma)}_{L^2(\RT)}\\
=&\rbr{w,\partial_{a_1} Z(\bm{a},\gamma)}_{L^2(\RT)}=\rbr{w,\partial_{a_2}Z(\bm{a},\gamma)}_{L^2(\RT)}=0,
\end{align*}
then 
\[\tbr{S_{4/\sqrt{3}}''(Z(\bm{a},\gamma))w,w}_{\mathcal{Z}^{-1},\mathcal{Z}^1} \geq k_1 \norm{w}_{\mathcal{Z}^1}^2,\]
where $\tbr{\cdot,\cdot}_{\mathcal{Z}^{-1},\mathcal{Z}^1}$ denotes the pairing between $\mathcal{Z}^{-1}$ and $\mathcal{Z}^1$ such that 
\[\tbr{w,v}_{\mathcal{Z}^{-1},\mathcal{Z}^{1}}=(w,v)_{L^2(\RT)}\]
for $w,v \in \mathcal{Z}^{1}$.
\end{lemma}
\proof
By the Fourier expansion in the transverse direction $y$, we have
\[S_{4/\sqrt{3}}''(Q_{4/\sqrt{3}})w=\sum_{n \in \Z} \Bigl(-\partial_x^2 w_n -n^{2}\partial_x^{-2}w_n + \frac{4}{\sqrt{3}} w_n - Q_{4/\sqrt{3}}w_n \Bigr)e^{iny},\]
where
\[w(x,y)=\sum_{n \in \Z} w_n(x)e^{iny}.\]
From Lemma 2.1 and Lemma 2.3 in \cite{R T 3}, there exists $C>0$ such that for $n \in \Z \setminus  \{-1,1\}$, $w_0 \in H^1(\R)$ and $w_n \in \mathcal{Z}_x^1$ with $(w_0,Q_{4/\sqrt{3}})_{L^2(\R)}=(w_0,\partial_xQ_{4/\sqrt{3}})_{L^2(\R)}=0$ we have
\begin{align}\label{eq-lem-5--1}
 \Tbr{\Bigl(-\partial_x^2+ \frac{4}{\sqrt{3}}-Q_{4/\sqrt{3}}\Bigr)w_0,w_0}_{H^{-1}(\R),H^1(\R)} \geq C \norm{w_0}_{H^1(\R)}^2
 \end{align}
and
\begin{align}\label{eq-lem-5-0}
 \Tbr{\Bigl(-\partial_x^2- n^2 \partial_x^{-2}+\frac{4}{\sqrt{3}}-Q_{4/\sqrt{3}}\Bigr)w_n,w_n}_{\mathcal{Z}^{-1}_x,\mathcal{Z}^1_x} \geq C(\norm{w_n}_{H^1(\R)}^2+n^2\norm{\partial_x^{-1}w_n}_{L^2(\R)}^2)
\end{align}
where 
\[\mathcal{Z}_x^s=\{u: \norm{(|\xi|+|\xi|^{-1})^s\hat{u}(\xi)}_{L^2_{\xi}(\R)}<\infty\}\]
and $\tbr{\cdot,\cdot}_{\mathcal{Z}^{-1}_x,\mathcal{Z}^1_x}$ denote the pairing between $\mathcal{Z}_x^{-1}$ and $\mathcal{Z}_x^1$.
Let $u_1=\partial_x^{-1}w_1$ and $\mathcal{L}_1=-\partial_x^2-\partial_x^{-2}+\frac{4}{\sqrt{3}}-Q_{4/\sqrt{3}}$.
Then, we have
\[ \Tbr{\mathcal{L}_1w_1,w_1 }_{\mathcal{Z}^{-1}_x,\mathcal{Z}^1_x} = \Tbr{-\partial_x\mathcal{L}_1\partial_xu_1,u_1}_{H^{-2}(\R),H^2(\R)}.\]
By applying the Weyl Lemma, we obtain that the essential spectrum of the operator $ -\partial_x\mathcal{L}_1\partial_x $ is $[n^2,\infty)$.
Thus, to show the coerciveness of $\mathcal{L}_1$ on a subspace, we investigate the eigenvalue problem
\begin{align}\label{eq-lem-5-1}
 -\partial_x\mathcal{L}_1\partial_x u_1=\lambda u_1, \quad u_1 \in H^4(\R).
\end{align}
From the argument in \cite{A P S,R T 3}, it was proved that the problem \eqref{eq-lem-5-1} has no negative eigenvalue.
Moreover, $-\partial_x \mathcal{L}_1 \partial_x u_n=0$ has linearly independent solutions 
\begin{align*}
\partial_xg_{\sqrt{3}}, \quad \partial_xg_{-\sqrt{3}}, \quad \partial_xg_1, \quad \lim_{\mu \to 1} \frac{\partial_x(g_\mu+g_{-\mu})}{\mu-1},
\end{align*}
where
\[g_{\mu}(x)=e^{3^{-1/4}\mu x}\bigl(\mu^3+2\mu-3\mu^2 \tanh(3^{-1/4}x)\bigr).\]
Then, 
\[\partial_x g_1(x)= -\frac{6 \sinh \frac{x}{3^{1/4}}}{3^{1/4}\cosh^2\frac{x}{3^{1/4}}}= \frac{1}{2\sqrt{2}}\partial_x^{-1}v_*. \]
Since the growth rates of $\partial_xg_{\sqrt{3}},   \partial_xg_{-\sqrt{3}} $ and $ \lim_{\mu \to 1} \frac{\partial_x(g_\mu+g_{-\mu})}{\mu-1}$ at the infinity are different and do not belong to $L^2$, we obtain 
\[L^2(\R) \cap \Bigl\{  c_1\partial_xg_{\sqrt{3}}+ c_2\partial_xg_{-\sqrt{3}} + c_3 \lim_{\mu \to 1} \frac{\partial_x(g_\mu+g_{-\mu})}{\mu-1}: c_1,c_2,c_3 \in \R \Bigr\}=\{0\}.\]
Therefore, the zero eigenvalue of 
%$-\partial_x^4-1+ \frac{4}{\sqrt{3}}\partial_x^2-\partial_x Q_{4/\sqrt{3}} \partial_x  $
$-\partial_x \mathcal{L}_1\partial_x$
 is simple and there exists $C>0$ such that 
\begin{align} \label{eq-lem-5-2}
\Tbr{\mathcal{L}_1w_1,w_1 }_{\mathcal{Z}^{-1}_x,\mathcal{Z}^1_x} = & \Tbr{-\partial_x \mathcal{L}_1\partial_xu_1,u_1}_{H^{-2}(\R),H^2(\R)} \notag \\
\geq & C' \norm{u_1}_{H^2}^2 \geq C (\norm{w_1 }_{H^1}^2+\norm{\partial_x^{-1}w_1}_{L^2}^2)
\end{align}
for $w_1 \in H^{1}(\R)$ with $\partial_x^{-1}w_1=u_1 \in L^2(\R)$ and $(w_1 ,v_*)_{L^2(\R)}=0$.
From \eqref{eq-lem-5--1}, \eqref{eq-lem-5-0} and \eqref{eq-lem-5-2}, we obtain that there exists $C_0>0$ such that 
\[\tbr{S_{4/\sqrt{3}}''(Q_{4/\sqrt{3}})w,w}_{\mathcal{Z}^{-1},\mathcal{Z}^1} \geq C_0 \norm{w}_{\mathcal{Z}^1}^2.\]
for $w \in \mathcal{Z}^1$ with
\[\rbr{w,Q_{4/\sqrt{3}}}_{L^2(\RT)}=\rbr{w,\partial_x Q_{4/\sqrt{3}}}_{L^2(\RT)}=\rbr{w,v_* \cos y}_{L^2(\RT)}=\rbr{w, v_* \sin y}_{L^2(\RT)}=0.\]
By a continuity argument and the above inequality, we obtain the conclusion.
\qed

In the following lemma, we set a modulation of the scaling in $x$ direction and show the expansion of the modulation.
\begin{lemma}\label{lem-1}
Let $B_1=\{ \bm{b} \in \R^2: |\bm{b}|<1\}$ and
\[\gamma_l(a)= \Bigl( M(Z(l)) M(Z(a))^{-1} \Bigr)^{\frac{2}{3}}\]
for $l\geq 0, a\in \R$.
For $\bm{a} \in B_1$
\[M(Z(\bm{a},\gamma_l(|\bm{a}|)))= M(Z(l)),\]
\begin{align}\label{eq-lem-1-1}
\gamma_{0}(|\bm{a}|)-1=-\frac{\partial_a^4M(Z(a))|_{a=0}}{36M(Q_{4/\sqrt{3}})}|\bm{a}|^4+o(|\bm{a}|^4) \mbox{ as } |\bm{a}|\to 0
\end{align}
and
\begin{align}\label{eq-lem-1-2}
\gamma_{l}(|\bm{a}|)-1=-\frac{2\partial_aM(Z(a))|_{a=l}}{3M(Z(l))}(|\bm{a}|-l)+O((|\bm{a}|-l)^2) \mbox{ as } |\bm{a}|\to l.
\end{align}
\end{lemma}
\proof
The equation \eqref{eq-lem-1-2} follows the Taylor expansion.
Since
\[M(Z(a))=M(Q_{4/\sqrt{3}})+\frac{\partial_a^4M(Z(a))|_{a=0}a^4}{4!}+o(a^4)\]
and
\begin{align*}
\gamma_0(|\bm{a}|)=&1-\frac{M(Z(|\bm{a}|))^{\frac{2}{3}}-M(Q_{4/\sqrt{3}})^{\frac{2}{3}}}{M(Z(|\bm{a}|))^{\frac{2}{3}}}=1-\frac{\partial_a^4M(Z(a))|_{a=0}}{36M(Q_{4/\sqrt{3}})}|\bm{a}|^4+o(|\bm{a}|^4),
\end{align*}
we have the equation \eqref{eq-lem-1-1}.

\qed

Next, we investigate the expansion of the Lyapunov function along the bifurcation parameter $a$. 
\begin{lemma}\label{lem-2}
Let $\epsilon _0>0$.
For $\bm{a} \in B_1 \setminus \{(0,0)\}$ and $0<l<\epsilon _0$,
\begin{align}\label{eq-lem-2-1}
S_{4/\sqrt{3}}\bigl(Z(\bm{a},\gamma_0(|\bm{a}|))\bigr)-S_{4/\sqrt{3}}(Q_{4/\sqrt{3}})=&\frac{5c''(0)\partial_a^4M(Z(a))|_{a=0}}{6!}|\bm{a}|^6 +o(|\bm{a}|^6)
\end{align}
and
\begin{align}\label{eq-lem-2-1-1}
&S_{c(l)}\bigl(Z(\bm{a},\gamma_l(|\bm{a}|))\bigr)-S_{c(l)}(Z(l)) \notag \\ 
=&\Bigl(\frac{(\partial_aM(Z(a))|_{a=l})^2\partial_c\norm{Q_c}_{L^2}^2|_{c=4/\sqrt{3}}}{9M(Z(l))^2}+\frac{c'(l)\partial_aM(Z(a))|_{a=l}}{2}\Bigr)(|\bm{a}|-l)^2 \notag \\
&+O((\epsilon _0 +||\bm{a}|-l|)(|\bm{a}|-l)^2).
\end{align}
\end{lemma}
\proof
By the expansion
\begin{align}\label{eq-lem-2-2}
Z(\bm{a},\gamma_{l}(|\bm{a}|))=&Z(\bm{a},1)+(\gamma_{l}(|\bm{a}|)-1)\partial_cQ_{4/\sqrt{3}}\notag \\
&+O\bigl((l+||\bm{a}|-l|+(\gamma_{l}(|\bm{a}|)-1))(\gamma_{l}(|\bm{a}|)-1)\bigr),
\end{align}
we have
\begin{align*}
&S_{c(l)}\bigl(Z(\bm{a},\gamma_{l}(|\bm{a}|))\bigr) -S_{c(l)}(Z(l))\\
=&S_{c(|\bm{a}|)}\bigl(Z(|\bm{a}|)\bigr) - S_{c(l)}(Z(l))+(c(l)-c(|\bm{a}|))M(Z(l))\\
&+\frac{1}{2}(\gamma_{l}(|\bm{a}|)-1)^2(S_{4/\sqrt{3}}''(Q_{4/\sqrt{3}})\partial_c Q_{4/\sqrt{3}}, \partial_c Q_{4/\sqrt{3}})_{L^2}+O((l+||\bm{a}|-l|)(\gamma_{l}(|\bm{a}|)-1)^2)
\end{align*}
for $l\geq 0$.
Since
\begin{align*}
 &\partial_a^n\Bigl(  S_{c(a)}\bigl(Z(a)\bigr) - S_{c(l)}(Z(l))+(c(l)-c(a))M(Z(l))\Bigr)\\
=&-c^{(n)}(a)M(Z(l))+ \sum_{j=0}^{n-1}{}_{n-1}C_j c^{(j+1)}(a)\partial_a^{n-1-j}M(Z(a)),
 \end{align*}
we obtain 
\begin{align}\label{eq-lem-2-3}   
& S_{c(|\bm{a}|)}\bigl(Z(|\bm{a}|)\bigr) - S_{4/\sqrt{3}}(Q_{4/\sqrt{3}})+(4/\sqrt{3}-c(|\bm{a}|))M(Q_{4/\sqrt{3}})\notag \\  
=&\frac{5}{6!}c''(0)\partial_a^4M(Z(a))|_{a=0}|\bm{a}|^6+o(|\bm{a}|^6)
\end{align}
and
\begin{align}\label{eq-lem-2-3-1}  
& S_{c(|\bm{a}|)}\bigl(Z(|\bm{a}|)\bigr) - S_{c(l)}(Z(l))+(c(l)-c(|\bm{a}|))M(Z(l)) \notag \\  
=&\frac{1}{2}c'(l)\partial_aM(Z(a))|_{a=l}(|\bm{a}|-l)^2+o((|\bm{a}|-l)^2)
\end{align}
for $l\geq 0$.
The equations \eqref{eq-lem-2-1} and \eqref{eq-lem-2-1-1} follow \eqref{eq-lem-2-3} and \eqref{eq-lem-2-3-1}.
% By combining the relation
% \begin{align*}
% &S_c\bigl(Z(\bm{a},\gamma_c(|\bm{a}|)\bigr) -S_c(Q_c) \\
% =& \Bigl( \frac{\sqrt{3}c}{4}\Bigr)^{\frac{5}{2}}\Bigl( S_{4/\sqrt{3}}\bigl(Z(\bm{a},\gamma_{4/\sqrt{3}}(|\bm{a}|)\bigr)-S_{4/\sqrt{3}}(Q_{4/\sqrt{3}})\Bigr)+\Bigl( 1-\frac{3c^2}{16}\Bigr)\norm{\partial_x^{-1}\partial_yZ(\bm{a},\gamma_c(|\bm{a}|))}_{L^2}^2
% \end{align*}
% and the equations \eqref{eq-lem-2-2}--\eqref{eq-lem-2-3}, we get
% \begin{align*}
%  &S_c\bigl(Z(\bm{a},\gamma_c(|\bm{a}|))\bigr)-S_c(Q_c)\\
%  =& \Bigl( \frac{\sqrt{3}c}{4}\Bigr)^{\frac{5}{2}}\frac{5c''(0)\partial_a^4M(Z(a))|_{a=0}}{6!}|\bm{a}|^6  +\Bigl( 1-\frac{3c^2}{16}\Bigr)\norm{\partial_x^{-1}\partial_yZ(\bm{a},\gamma_c(|\bm{a}|))}_{L^2}^2+o(|\bm{a}|^6).
%  \end{align*}

\qed

% We define a distance from the Zaitsev solitons $\{Z(l)(\cdot-x_0,\cdot-y_0)\}_{x_0,y_0}$ by
Let a distance
\[\mbox{dist}_l(u)=\inf_{(x_0,y_0) \in \RT}\norm{u(\cdot,\cdot)-Z(l)(\cdot-x_0,\cdot-y_0)}_{\mathcal{Z}^1}\]
and neighborhoods 
\[N_{\varepsilon ,l}=\{u \in \mathcal{Z}^1: \mbox{dist}_l(u)<\varepsilon \},\]
\[N_{\varepsilon ,l}^k=\{u \in N_{\varepsilon ,l}: M(u)=M(Z(k))\}.\]
In the following lemma, we define the modulation parameters to control the degenerate direction of the Lyapunov function.
The following lemma follows the implicit function theorem.
\begin{lemma}\label{lem-3}
Let $\varepsilon >0$ sufficiently small.
Then, there exist $K_1>0$, $C^2$ functions $\rho:N_{\varepsilon ,0} \to \R$, $\gamma: N_{\varepsilon ,0} \to \R$, $\bm{a}=(a_1,a_2):N_{\varepsilon ,0} \to B_1$ and $\eta:N_{\varepsilon ,0} \to \mathcal{Z}^1$ such that for $u \in N_{\varepsilon ,0} $
\[u(\cdot+\rho(u),\cdot)=Z(\bm{a}(u),\gamma(u))(\cdot,\cdot)+\eta(u)(\cdot,\cdot),\]
\begin{align}\label{eq-lem-3-1}
|\gamma(u)-1|+|\bm{a}(u)|+\norm{\eta(u)}_{\mathcal{Z}^1}\leq K_1 \mbox{\rm dist}_{0}(u)
\end{align}
and
$\rbr{\eta(u),Z(\bm{a}(u),\gamma(u))}_{L^2}=\rbr{\eta(u),\partial_xZ(\bm{a}(u),\gamma(u))}_{L^2}=\rbr{\eta(u),\partial_{a_1}Z(\bm{a}(u),\gamma(u))}_{L^2}=\rbr{\eta(u),\partial_{a_2}Z(\bm{a}(u),\gamma(u))}_{L^2}=0.$
Moreover, $\gamma, \rho$ and $ \bm{a}=(a_1,a_2)$ satisfy
\[G(u,\gamma(u),\rho(u),a_1(u),a_2(u))=0,\]
where 
\[G(u,\gamma,\rho,a_1,a_2)=
\begin{pmatrix}
\rbr{u(\cdot-\rho,\cdot)-Z(\bm{a},\gamma),Z(\bm{a},\gamma)}_{L^2} \\
\rbr{u(\cdot-\rho,\cdot)-Z(\bm{a},\gamma),\partial_xZ(\bm{a},\gamma)}_{L^2}\\
\rbr{u(\cdot-\rho,\cdot)-Z(\bm{a},\gamma),\partial_{a_1} Z(\bm{a},\gamma)}_{L^2}\\
\rbr{u(\cdot-\rho,\cdot)-Z(\bm{a},\gamma),\partial_{a_2} Z(\bm{a},\gamma)}_{L^2}
\end{pmatrix}.\]
\end{lemma}
\proof
Since $G(1,0,0,0)=0$ and the Jacobian matrix
\[\frac{\partial G}{\partial(\gamma,\rho,a_1,a_2)}\Bigl|_{\gamma=1,\rho=0,a_1=0,a_2=0}\]
is regular, there exists $C^2$ functions $\gamma,\rho,a_1,a_2: N_{\varepsilon ,0} \to \R$ such that \[G(u,\gamma(u),\rho(u),a_1(u),a_2(u))=0.\]
Let 
\[\eta(u)(\cdot,\cdot)=u(\cdot+\rho(u),\cdot)-Z(\bm{a}(u),\gamma(u))(\cdot,\cdot).\]
Then, $\eta$ satisfies the orthogonal condition and \eqref{eq-lem-3-1}.
\qed

The following lemma shows the smallness of the difference between $Z(\bm{a}(u),\gamma_l(|\bm{a}(u)|))$ and $Z(\bm{a}(u),\gamma(u))$.
\begin{lemma}\label{lem-4}
Let $\varepsilon >0$ sufficiently small.
There exists $C>0$ such that for $0\leq l<2\varepsilon $ and $u \in N_{\varepsilon ,0}^l$
\[\norm{Z(\bm{a}(u),\gamma_l(|\bm{a}(u)|))-Z(\bm{a}(u),\gamma(u))}_{\mathcal{Z}^1}\leq C \norm{\eta(u)}_{L^2}^2,\]
\begin{align}\label{eq-lem-4-1}
|\gamma_l(\bm{a}(u))-\gamma(u)|\leq C \norm{\eta(u)}_{L^2}^2.
\end{align}
\end{lemma}
\proof
By Lemma \ref{lem-3}, 
\begin{align*}
M\bigl(Z(\bm{a}(u),\gamma_l(\bm{a}(u)))\bigr)=M(Z(l))=&M\bigl(\eta(u)+Z(\bm{a}(u),\gamma(u))\bigr)\\
=&M(\eta(u))+M\bigl( Z(\bm{a}(u),\gamma(u)) \bigr)
\end{align*}
for $u \in N_{\varepsilon ,0}^l$.
Since there exists $K>0$ such that 
\[ |\gamma(u)-1|+|\gamma_l(\bm{a}(u))-1|<K \varepsilon ,\]
we obtain
\begin{align*}
M(\eta(u))=M\bigl(Z(\bm{a}(u),\gamma_l(\bm{a}(u)))\bigr)-M\bigl(Z(\bm{a}(u),\gamma(u))\bigr)=&(\gamma_l(\bm{a}(u))^{\frac{3}{2}}-\gamma(u)^{\frac{3}{2}})M(Z(\bm{a}(u)))\\
\gtrsim & \gamma_l(\bm{a}(u))-\gamma(u)\geq 0.
\end{align*}
Therefore, we have
\[\norm{Z(\bm{a}(u),\gamma_l(\bm{a}(u)))-Z(\bm{a}(u),\gamma(u))}_{\mathcal{Z}^1} \lesssim \gamma_l(\bm{a}(u))-\gamma(u) \lesssim M(\eta(u)).\]
\qed

The following lemma shows that the degeneracy of the linearized operator $S_{4/\sqrt{3}}''(Q_{4/\sqrt{3}})$ is covered by the sixth order nonlinearity of $\bm{a}$.
\begin{lemma}\label{lem-6}
There exists $\varepsilon _1>0$ such that for $0<\varepsilon <\varepsilon _1$ and $u \in N_{\varepsilon ,0}^{0}$ 
\begin{align}\label{eq-lem-6-1}
S_{4/\sqrt{3}}(u)-S_{4/\sqrt{3}}(Q_{4/\sqrt{3}})  =&\frac{5c''(0)\partial_a^4M(Z(a))|_{a=0}}{6!}|\bm{a}(u)|^6 \notag \\
&+\frac{1}{2}\tbr{S_{4/\sqrt{3}}''\bigl(Z(\bm{a}(u),\gamma_{0}(|\bm{a}(u)|))\bigr)\eta(u),\eta(u)}_{\mathcal{Z}^{-1},\mathcal{Z}^1}\notag \\
&+o(\norm{\eta(u)}_{\mathcal{Z}^1}^2+|\bm{a}(u)|^6).
\end{align}
Moreover, for $0<\varepsilon <\varepsilon _1$, $l>0$ and $u \in N_{\varepsilon ,0}^{l}$ 
\begin{align}\label{eq-lem-6-2}
&S_{c(l)}(u)-S_{c(l)}(Z(l)) \notag \\
&= \Bigl(\frac{(\partial_aM(Z(a))|_{a=l})^2\partial_c\norm{Q_c}_{L^2}^2|_{c=4/\sqrt{3}}}{9M(Z(l))^2}+\frac{c'(l)\partial_aM(Z(a))|_{a=l}}{2}\Bigr)(|\bm{a}|-l)^2\notag \\
&+\frac{1}{2}\tbr{S_{4/\sqrt{3}}''\bigl(Z(\bm{a}(u),\gamma_{0}(|\bm{a}(u)|))\bigr)\eta(u),\eta(u)}_{\mathcal{Z}^{-1},\mathcal{Z}^1}\notag \\
&+O((l+\norm{\eta(u)}_{\mathcal{Z}^1}+||\bm{a}(u)|-l|)(\norm{\eta(u)}_{\mathcal{Z}^1}^2+(|\bm{a}(u)|-l)^2)).
\end{align}
\end{lemma}
\proof 
From the Taylor expansion, we have for $l\geq 0$
\begin{align*}
&S_{c(l)}(u)-S_{c(l)}(Z(l)) \\
 =& S_{c(l)}(Z(\bm{a}(u),\gamma(u))+\eta(u)) -S_{c(l)}(Z(l))  \\
=&S_{c(l)}\bigl(Z(\bm{a}(u),\gamma_{l}(|\bm{a}(u)|))\bigr)-S_{c(l)}(Z(l))  \\
&+ \tbr{S_{c(l)}'\bigl(Z(\bm{a}(u),\gamma_{l}(|\bm{a}(u)|))\bigr),Z(\bm{a}(u),\gamma(u))+\eta(u) -Z(\bm{a}(u),\gamma_{l}(|\bm{a}(u)|)) }_{\mathcal{Z}^{-1},\mathcal{Z}^1}\\
&+\frac{1}{2}\tbr{S_{c(l)}''\bigl(Z(\bm{a}(u),\gamma_{l}(|\bm{a}(u)|))\bigr)\eta(u),\eta(u)}_{\mathcal{Z}^{-1},\mathcal{Z}^1}+o(\norm{\eta(u)}_{\mathcal{Z}^1}^2)
% \\
% =&\frac{5c''(0)\partial_a^4M(Z(a))|_{a=0}}{6!}|\bm{a}(u)|^6+\tbr{S_{4/\sqrt{3}}'\bigl(Z(\bm{a}(u),\gamma_{0}(|\bm{a}(u)|))\bigr),\eta(u)}_{\mathcal{Z}^{-1},\mathcal{Z}^1}\\
% &+\frac{1}{2}\tbr{S_{4/\sqrt{3}}''\bigl(Z(\bm{a}(u),\gamma_{0}(|\bm{a}(u)|))\bigr)\eta(u),\eta(u)}_{\mathcal{Z}^{-1},\mathcal{Z}^1}+o(\norm{\eta(u)}_{\mathcal{Z}^1}^2+|\bm{a}(u)|^6)
.
\end{align*}
Since $S_{4/\sqrt{3}}''(Q_{4/\sqrt{3}})\partial_cQ_{4/\sqrt{3}}=-Q_{4/\sqrt{3}}$, by \eqref{eq-lem-4-1} we have
\begin{align*}
&\tbr{S_{c(l)}'(Z(\bm{a}(u),\gamma_{l}(|\bm{a}(u)|)),\eta(u)}_{\mathcal{Z}^{-1},\mathcal{Z}^1}\\
% =&\tbr{S_{c(|\bm{a}(u)|)}'(Z(\bm{a}(u),\gamma_{l}(|\bm{a}(u)|)),\eta(u)}_{\mathcal{Z}^{-1},\mathcal{Z}^1}+(c(l)-c(|\bm{a}(u)|))(Z(\bm{a}(u),\gamma_{l}(|\bm{a}(u)|)),\eta(u))_{L^2}\\
=& \tbr{\bigl( S_{c(|\bm{a}(u)|)}''(Z(\bm{a}(u),\gamma(u)))-S_{4/\sqrt{3}}''(Q_{4/\sqrt{3}}) \bigr)(\gamma_{l}(\bm{a}(u))-\gamma(u))\partial_cQ_{4/\sqrt{3}},\eta(u)}_{\mathcal{Z}^{-1},\mathcal{Z}^1}\\
&- (\gamma_{l}(|\bm{a}(u)|)-\gamma(u))(Q_{4/\sqrt{3}},\eta(u))_{L^2}\\
&+(c(l)-c(|\bm{a}(u)|))\bigl(Z(\bm{a}(u),\gamma_{l}(|\bm{a}(u)|))-Z(\bm{a}(u),\gamma(u)),\eta(u)\bigr)_{L^2}\\
&+O\bigl((l+||\bm{a}|-l|+\norm{\eta(u)}_{\mathcal{Z}^1})(|\gamma_l(|\bm{a}(u)|)-\gamma(u)|^2+\norm{\eta(u)}_{\mathcal{Z}^1}^2)\bigr)\\
=&O\bigl((l+||\bm{a}|-l|+\norm{\eta(u)}_{\mathcal{Z}^1})(|\gamma_l(|\bm{a}(u)|)-\gamma(u)|^2+\norm{\eta(u)}_{\mathcal{Z}^1}^2)\bigr).
\end{align*}
Therefore, by Lemma \ref{lem-2} we obtain \eqref{eq-lem-6-1} and \eqref{eq-lem-6-2}.
\qed

\section{Proof of the orbital stability of the line soliton with the critical speed}
In this section, applying Lemma \ref{lem-5} and Lemma \ref{lem-6}, we show the orbital stability of the line soliton $Q_{4/\sqrt{3}}$ with the critical traveling speed.
We prove the orbital stability by the contradiction.
We assume the line soliton $Q_{4/\sqrt{3}}$ is unstable.
Then, there exist $\varepsilon _0>0$,  a sequence $\{u_n\}_n$ of solutions to \eqref{KP-I-eq} and a sequence $\{t_n\}_n$ such that $t_n>0$, $u_n(0) \to Q_{4/\sqrt{3}}$  as $n \to \infty$ in $\mathcal{Z}^1$ and $\mbox{dist}_{0}(u_n(t_n))>\varepsilon _0$.
Let $v_n = M(Q_{4/\sqrt{3}})^{1/2}M(u_n)^{-1/2}u_n(t_n)$.
Then, we have $M(v_n)=M(Q_{4/\sqrt{3}})$, $\lim_{n \to \infty}\norm{v_n -u_n(t_n)}_{\mathcal{Z}^{1}}=0$ and $\lim_{n \to \infty}S_{4/\sqrt{3}}(v_n)=S_{4/\sqrt{3}}(Q_{4/\sqrt{3}})$.
By Lemma \ref{lem-5} and Lemma \ref{lem-6}, we obtain that there exists $k_0>0$ such that
\begin{align}\label{eq-sta-1}
S_{4/\sqrt{3}}(v_n)-S_{4/\sqrt{3}}(Q_{4/\sqrt{3}}) \geq k_0 (|\bm{a}(v_n)|^6+\norm{\eta(v_n)}_{\mathcal{Z}^1}^2)
\end{align}
for sufficiently large $n$.
From the inequality \eqref{eq-sta-1}, we have
\begin{align}\label{eq-sta-2}
\lim_{n \to \infty} \bm{a}(v_n)=0, \quad \lim_{n \to \infty}\norm{\eta(v_n)}_{\mathcal{Z}^1}=0.
\end{align}
By Lemma \ref{lem-4}, the equation $\lim_{n \to \infty} \gamma_{0}(\bm{a}(v_n))=1$ yields 
\begin{align}\label{eq-sta-3}
\lim_{n \to \infty} \gamma(v_n) =1.
\end{align}
From \eqref{eq-sta-2} and \eqref{eq-sta-3}, we obtain 
\[\lim_{n \to \infty} \mbox{dist}_{0}(u_n(t_n)) =\lim_{n \to \infty} \mbox{dist}_{0}(v_n)=0.\]
This is a contradiction.
Thus, $Q_{4/\sqrt{3}}$ is orbitally stable.

\section{Proof of the orbital stability of the Zaitsev solitons near by the line soliton}
In this section, we show the orbital stability of the Zaitsev solitons near by the line soliton $Q_{4/\sqrt{3}}$.
The positivity of $c'(l)$ and $\partial_aM(Z(l))$ for sufficiently small $l>0$ follows
 $c''(0)=4/\sqrt{3}$ and $\partial_a^4M(Z(a))|_{a=0}>0$.
Thus, from Lemma \ref{lem-5} and Lemma \ref{lem-6}, we have that for sufficiently small $l>0$  there exist $k_l,\varepsilon _0>0$ such that
\[S_{c(l)}(u)-S_{c(l)}(Z(l)) \geq k_l((|\bm{a}(u)|-l)^2+\norm{\eta(u)}_{\mathcal{Z}^1}^2)\]
for $u \in N_{\varepsilon_0,l}^l$.
Therefore, applying the same argument as the argument in the proof of the orbital stability of $Q_{4/\sqrt{3}}$, we obtain the orbital stability of the Zaitsev soliton $Z(l)$ for sufficiently small $l>0$.

\vspace{0.5cm}

% Funding: 
% This research did not receive any specific grant from funding agencies in the public, commercial, or
% not-for-profit sectors. 

\section*{Acknowledgments}
The author is grateful to Nikolay Tzvetkov for his helpful advices.

\end{document}